\documentclass[12pt,twoside]{article}
\usepackage{amsmath, amssymb, amsthm}
\newtheorem{lemma}{Lemma}

\font\fourteenb=cmb10 at 14pt
\begin{document}
\vspace*{-1.0cm}\noindent \copyright
 Journal of Technical University at Plovdiv\\[-0.0mm]\
\ Fundamental Sciences and Applications, Vol. 9, 2000\\[-0.0mm]
\textit{Series A-Pure and Applied Mathematics}\\[-0.0mm]
\ Bulgaria, ISSN 1310-8271\\[+1.2cm]
\font\fourteenb=cmb10 at 14pt
\begin{center}

   {\bf \LARGE Composition operators on the space of
Couchy-Stiltjes transforms \\ \ \\ \large Peyo Stoilov and
Roumyana Gesheva }
\end{center}

\

\


\begin{abstract}
Let $K\,$ denote the space of all Couchy-Stiltjes transforms. Let
$\varphi\,$ be an analytic map of the unit disk into itself and
$C_{\varphi}\,$ denote the composition operator in $K$. In this
note is given a new proof of the norm estimate of J. Cima and A.
Matheson:
$$
||C_{\varphi}||_{K}\leq\frac{1+2|\varphi(0)|}{1-|\varphi(0)|}.
$$
\end{abstract}

\section{Introduction}
\hspace{0.2in} Let $K\,$ denote the family of all functions $f$,
analytic in the unit disk $\Bbb D\,$ and for which there exists a
finite Borel measure $\mu\,$ on the unit circle $\Bbb T\,$ such
that
$$
f(z)=\int\limits_{\Bbb
T}\frac{d\mu(\zeta)}{1-\overline{\zeta}z}\,\stackrel{\rm
def}{=}\,K_{\mu}.
$$
$K\,$ is a Banach space with the natural norm
$$
||f||_{K}=\inf\{||\mu||:\mu\in M, K_{\mu}=f\}.
$$
Let $\varphi\,$ be an analytic map of the unit disk into itself
and $C_{\varphi}\,$ denote the composition operator in $K$:
$$
C_{\varphi}f=f\circ\varphi,\quad f\in K.
$$
In [1] P. Bourdon and J. Cima proved that $C_{\varphi}\,$ is a
bounded operator on $K\,$ and
$$
||C_{\varphi}||_{K}\leq\frac{2+2\sqrt{2}}{1-|\varphi(0)|}.
$$
In [2] J. Cima and A. Matheson proved the norm estimate
\begin{equation}
||C_{\varphi}||_{K}\leq\frac{1+2|\varphi(0)|}{1-|\varphi(0)|}.
\end{equation}
This estimate is sharp in th sense that there are functions
$\varphi\,$ with $\varphi(0)\neq0\,$ for which equality is
attained.

The proof of Cima and Matheson is based on the following two
lemmas.
\begin{lemma}
If $f\in K, \psi\,$ is analytic in $\Bbb{D}\,$ and
$\psi(0)=0,\,|\psi(z)|<1\,$ in $\Bbb{D}$, then $f\circ\psi\in K\,$
and
$$
||f\circ\psi||_K\leq||f||_K.
$$
\end{lemma}
\begin{lemma}
For every $a\in \Bbb{D}$, let
$\lambda_{a}(z)=\frac{a-z}{1-\overline{a}z}$. Then
$f\circ\lambda_{a}\in K\,$ for every $f\in K\,$ and
\begin{equation}
||f\circ\lambda_{a}||_K\leq\frac{1+2|a|}{1-|a|}||f||_K.
\end{equation}
\end{lemma}
Lemma 2 leads to a quick proof of the estimate (1). Indeed, if
$\varphi\,$ is analytic map of $\Bbb{D}\,$ into itself, then
$$
\psi(z)=\frac{\varphi(0)-\varphi(z)}{1-\overline{\varphi(0)}\varphi(z)}
$$
is analytic map of $\Bbb{D}\,$ into itself and $\psi(0)=0$. Since
$$
\varphi(z)=\frac{\varphi(0)-\psi(z)}{1-\overline{\varphi(0)}\psi(z)}=\lambda_{a}\circ\psi(z),\quad
a=\varphi(0),
$$
then
$$
||f\circ\varphi||_{K}=||(f\circ\lambda_{a})\circ\psi(z)||_{K}\leq
||f\circ\lambda_{a}||_K\leq\frac{1+2|a|}{1-|a|}||f||_K
$$
for every $f\in K$. Hence
$$
||C_{\varphi}||_{K}\leq
\frac{1+2|a|}{1-|a|}=\frac{1+2|\varphi(0)|}{1-|\varphi(0)|}.
$$
The motivation for this paper is the following new proof of Lemma
2.
\section{New proof of Lemma 2.}
\hspace{0.2in} Let $C_{A}\,$ denote the space of all functions
analytic in $\Bbb{D}\,$ and continuous on $\overline{\Bbb{D}}\,$
for which $||f||_{C_A}=||f||_{\infty}$. Let
$\lambda(z)=\lambda_{a}(z),\; a\in\Bbb{D}$.\newline The estimate
(2) will be proved by using that $K\cong C^{*}_{A}\,$ under the
paring
$$
\langle f,h\rangle=\lim\limits_{r\to
1}\int\limits_{\Bbb{T}}f(rt)\overline{h(t)}\,dm(t),
$$
where $f\in K,\; h\in C_{A}\,$ and $dm(t)=\frac{1}{2\pi
i}\frac{dt}{t}$.\newline If $f=K_{\mu}\in K\,$ and $h\in C_{A}$,
then
\begin{eqnarray*}
\int\limits_{\Bbb{T}}(f\circ\lambda)\overline{h(t)}\,dm(t) &=&
\int\limits_{\Bbb{T}}\left(\int\limits_{\Bbb{T}}\frac{d\mu(\zeta)}
{1-\overline{\zeta}\lambda(rt)}\right)\overline{h(t)}\,dm(t)=\\
 & =& \int\limits_{\Bbb{T}}\left(\int\limits_{\Bbb{T}}\frac{\overline{h(t)}}
{1-\overline{\zeta}\lambda(rt)}dm(t)\right)\,d\mu(\zeta)
\end{eqnarray*}
and it follows that
$$
\left|\int\limits_{\Bbb{T}}(f\circ\lambda)\overline{h(t)}\,dm(t)\right|\leq
|\mu|
 \left\|\int\limits_{\Bbb{T}}\frac{h(t)}
{1-\zeta\overline{\lambda(rt)}}dm(t)\right\|_{\infty},
$$
which implies
\begin{equation}
||f\circ\lambda||_{K}\leq||f||_{K}\sup\left\{\lim\limits_{r\to
1}\left\|\int\limits_{\Bbb{T}}\frac{h(t)}
{1-\zeta\overline{\lambda(rt)}}dm(t)\right\|_{\infty}\;:\;||h||_{\infty}\leq
1\right\}.
\end{equation}
Since
$$
\frac{1}{1-\zeta\overline{\lambda(rt)}}=\frac{t-ar}{(1-\zeta\overline{a})t+r(\zeta-a)},
$$
applying Cauchy's theorem, we have
\begin{eqnarray*}
\int\limits_{\Bbb{T}}\frac{h(t)}{1-\zeta\overline{\lambda(rt)}}dm(t)
&=&
\frac{1}{2\pi i}\int\limits_{\Bbb{T}}h(t)\frac{t-ar}{(1-\zeta\overline{a})t+r(\zeta-a)}\,\frac{dt}{t}=\\
    & =& -\frac{arh(0)}{r(\zeta-a)}+h\left(r\frac{a-\zeta}{1-\zeta\overline{a}}\right)
\frac{r\frac{a-\zeta}{1-\zeta\overline{a}}-ar}{r(a-\zeta)}=\\
    & =& -a\frac{h(0)}{\zeta-a}+h\left(r\frac{a-\zeta}{1-\zeta\overline{a}}\right)\frac{1-|a|^2}{|1-\zeta\overline{a}|^2}.
\end{eqnarray*}
If $||h||_{\infty}\leq 1$, then
$$
 \left|\int\limits_{\Bbb{T}}\frac{h(t)}
{1-\zeta\overline{\lambda(rt)}}dm(t)\right|\leq
\frac{|a|}{1-|a|}+\frac{1-|a|^2}{(1-|a|)^2}=\frac{1+2|a|}{1-|a|}
$$
and from inequality (3) we obtain
$$
||f\circ\lambda||_{K}\leq \frac{1+2|a|}{1-|a|}||f\||_{K}.
$$
\section{Remarks}
\hspace{0.2in} Let $\varphi\,$ be an analytic map of the unit disk
$\Bbb{D}\,$ into itself and define an operator $P_{\varphi}\,$ on
$C_A\,$ by
$$
P_{\varphi}h=\lim\limits_{r\to 1}\int\limits_{\Bbb{T}}\frac{h(t)}
{1-\zeta\overline{\varphi(rt)}}dm(t)\,,\quad h\in C_A.
$$
From the proof of Lemma 2 and inequality (1) it follows that
$P_{\varphi}\,$ is a bounded operator on $C_A\,$ and
$$
||P_{\varphi}||_{C_{A}}\leq
\frac{1+2|\varphi(0)|}{1-|\varphi(0)|}.
$$
We shall note that $P_{\varphi}\,$ is an integral analogue of the
operator $\tilde{A}_{\varphi}$, applied in
[2].\vspace{2cm}\newline {\bf REFERENCES}\vspace{0.5cm}\newline 1.
P. Bourdon and J. Cima, {\it On integrals of Couchy--Stiltjes
type}, Houston J. Math. 14 (1988), 465-474.\newline 2. J. Cima and
A. Matheson, {\it Couchy transform and Composition operators},
Illinois J. Math. 42 (1998), 58-69.\vspace{0.5cm}\newline
\noindent
{\small Department of Mathematics\\
Technical University\\
25, Tsanko Dijstabanov,\\
Plovdiv, Bulgaria\\
e-mail: peyyyo@mail.bg}

\end{document}